\documentclass[letterpaper,11pt]{amsart}

\usepackage{latexsym,array,delarray,amsthm,amssymb,epsfig}

\theoremstyle{plain}
\newtheorem{theorem}{Theorem}[section]
\newtheorem{lemma}[theorem]{Lemma}
\newtheorem{proposition}[theorem]{Proposition}
\newtheorem{corollary}[theorem]{Corollary}

\theoremstyle{definition}
\newtheorem{definition}[theorem]{Definition}

\newtheorem{example}[theorem]{Example}
\newtheorem{problem}[theorem]{Problem}
\newtheorem{remark}[theorem]{Remark}

\newtheorem{algorithm}[theorem]{Algorithm}

\newcommand{\zz}{\mathbb{Z}}

\newcommand{\pp}{\mathbb{P}}

\newcommand{\qq}{\mathbb{Q}}
\newcommand{\rr}{\mathbb{R}}
\newcommand{\cc}{\mathbb{C}}

\newcommand{\xx}{\mathbf{x}}

\newcommand{\TT}{\mathcal{T}}

\def\trop{\ensuremath{\mathcal{T}}}
\def\myinit{\textup{in}}
\def\F{\ensuremath{\mathcal{F}}}
\def\G{\ensuremath{\mathcal{G}}}
\def\v{\ensuremath{{\bf{v}}}}
\def\u{\ensuremath{{\bf{u}}}}
\def\idim{\textup{dim}}
\def\ihom{\textup{homog}}
\def\andersomega{{w}}
\def\myinit{\textup{in}}
\def\TAB{\hspace*{0.5cm}}

\begin{document}
\title{Computing Tropical Varieties}

\author{
T. Bogart,\
A. Jensen,\
D. Speyer,\
B. Sturmfels,\
R. Thomas
}

\begin{abstract}
  The tropical variety of a $d$-dimensional prime ideal in a
  polynomial ring with complex coefficients is a pure $d$-dimensional
  polyhedral fan. This fan is shown to be connected in codimension
  one. We present algorithmic tools for computing the tropical
  variety, and we discuss our implementation of these tools in the
  Gr\"obner fan software \texttt{Gfan}.  Every ideal is shown to have
  a finite tropical basis, and a sharp lower bound is given for the
  size of a tropical basis for an ideal of linear forms.
\end{abstract}
\maketitle


\section{Introduction}

Every ideal in a polynomial ring with complex coefficients defines a tropical variety, which is
a polyhedral fan in a real vector space.  The objective of this
 paper is to introduce methods for computing
this fan, which coincides with the
``logarithmic limit set'' in George Bergman's seminal paper \cite{Ber}.

Given any polynomial $f \in \cc[x_1,x_2,\ldots,x_n]$ and a vector $w
\in \rr^n$, the {\rm initial form} ${\rm in}_w(f)$ is the sum of all
terms in $f$ of lowest $w$-weight; for instance, if $\, \ell = x_1 +
x_2 + x_3 + 1\,$ then $\,{\rm in}_{(0,0,1)} (\ell) \, = x_1+x_2+1\,$
and $\, {\rm in}_{(0,0,-1)}(\ell) \, = \, x_3 $.  The {\em tropical
  hypersurface} of $f$ is the set
$$\, \mathcal{T}(f) \quad = \,\quad \{ \, w \in \rr^n \,: \, {\rm
  in}_w(f)\,\,\, \hbox{is not a monomial}\,\} . $$
Equivalently,
$\mathcal{T}(f)$ is the union of all codimension one cones in the
inner normal fan of the Newton polytope of $f$.  Note that
$\mathcal{T}(f)$ is invariant under dilation, so we may specify
$\mathcal{T}(f)$ by giving its intersection with the unit sphere.  For
the linear polynomial $\ell$ above, $\mathcal{T}(\ell)$ is a
two-dimensional fan with six maximal cones. Its intersection with the
$2$-sphere is the complete graph on the four nodes $(1,0,0)$,
$(0,1,0)$, $(0,0,1)$ and
$-(\frac{1}{\sqrt{3}},\frac{1}{\sqrt{3}},\frac{1}{\sqrt{3}})$.

A finite intersection of tropical hypersurfaces is a {\em
  tropical prevariety} \cite{RST}. If we pick the second linear
form $\ell' = x_1 + x_2 + 2 x_3$ then 
$\mathcal{T}(\ell')$ is a graph with two vertices connected by three
edges on the $2$-sphere, and $\,\mathcal{T}(\ell) \cap
\mathcal{T}(\ell')$ consists of three edges of $\mathcal{T}(\ell)$
which are adjacent to
$-(\frac{1}{\sqrt{3}},\frac{1}{\sqrt{3}},\frac{1}{\sqrt{3}})$.  In
particular, the tropical prevariety $\,\mathcal{T}(\ell) \cap
\mathcal{T}(\ell')$ is not a tropical variety.

Tropical varieties are derived from ideals.  Namely, if $I$ is an
ideal in $\cc[x_1,\ldots,x_n]$ then its {\em tropical variety}
$\TT(I)$ is the intersection of the tropical hypersurfaces
$\mathcal{T}(f)$ where $f$ runs over all polynomials in $I$.
  Theorem~\ref{thm:tropical bases} below states that every
tropical variety is actually a tropical prevariety, i.e., the ideal
$I$ has a finite generating set $ \{f_1,f_2,\ldots,f_r \} $ such that
$$ \mathcal{T}(I) \quad = \quad
\mathcal{T}(f_1) \,\cap \,
\mathcal{T}(f_2) \, \cap \,\cdots \, \cap
\mathcal{T}(f_r) \, .$$
If this holds then $\{f_1,f_2,\ldots,f_r \}$ is called
a {\em tropical basis} of $I$. For
instance, our ideal $I = \langle \ell, \ell' \rangle $ 
 has the tropical basis $\,\{\, x_1 + x_2 + 2 x_3,\,
x_1+x_2 + 2, \,x_3 - 1\, \}$, and we find that its tropical variety
consists of three points on the sphere: 
$$ \mathcal{T}(I) \,\,\,\, = \,\,\,
\bigl\{ (1,0,0) ,\,(0,1,0),\, -(\frac{1}{\sqrt{2}},\frac{1}{\sqrt{2}},0) \bigr\}.$$

Our main contribution is a practical algorithm, along with its
implementation, for computing the tropical variety $\mathcal{T}(I)$
from any generating set of its ideal $I$. The emphasis lies on the geometric
and algebraic features of this computation. We do not address
issues of computational complexity, which have been studied by
Theobald \cite{The}. Our paper is organized as follows.

In Section 2 we give precise specifications of the algorithmic
problems we are dealing with, including the computation of a tropical
basis. We show that a finite tropical basis exists for every ideal
$I$, and we give tight bounds on its size for
linear ideals, thereby answering the question raised in \cite[\S 5,
page 13]{SS1}.  In Section 3 we prove that the tropical variety
$\TT(I)$ of a prime ideal $I$ is connected in codimension one. This result
is the foundation of Algorithm~\ref{alg: traverse} for computing
$\TT(I)$.  Section 4 also describes methods for computing tropical
bases and tropical prevarieties. Our algorithms have been
implemented in the software package {\tt Gfan} \cite{Gfan}. In Section
5 we compute the tropical variety of several non-trivial ideals using
{\tt Gfan}. The tropical variety $\TT(I)$ is a subfan of the Gr\"obner
fan of $I$ (defined in Section~2). The Gr\"obner fan is generally much
more complicated and harder to compute than $\TT(I)$. In Section 6 we
compare these two fans, and we exhibit a family of curves for which
the tropical variety of each member consists of four rays but the
number of one-dimensional cones in the Gr\"obner fan grows
arbitrarily.

A note on the choice of ground field is in order. In this paper
we will work with varieties defined over $\cc$. In the implementation
of our algorithm (Section \ref{sec:examples}), we have required our
polynomials to have rational coefficients,
 but our algorithms do not use any particular properties
of $\qq$. It is important, however, that we work over a field of
characteristic $0$, as our proof of correctness uses the
Kleiman-Bertini theorem in the proof of Theorem   \ref{codimone}.

In most papers on tropical algebraic geometry (cf.~\cite{EKL, Mik,
  RST, SS, The}), tropical varieties are defined from polynomials with
coefficients in a field $K$ with a non-archimedean valuation. These
tropical varieties are not fans but polyhedral complexes.  We close
the introduction by illustrating how our algorithms can be applied to
this situation. Consider the field $ \cc(\epsilon)$ of rational
functions in the unknown $\epsilon$.  Then $\cc(\epsilon)$ is a
subfield of the algebraically closed field $\cc \{ \! \{ \epsilon \}
\! \}$ of Puiseux series with real exponents, which is an example of a
field $K$ as in the above cited papers.  Suppose we are given an ideal
$I$ in $\cc(\epsilon)[x_1, \ldots,x_n]$. Let $I' \subset \cc \{ \! \{
\epsilon \} \! \}[x_1, \ldots, x_n]$ be the ideal generated by $I$.
The tropical variety $\mathcal{T}(I')$, in the sense of the papers
above, is a finite polyhedral complex in $\rr^n$ which usually has both
bounded and unbounded faces. To study this complex, we consider the
polynomial ring in $n+1$ variables, $\,\cc[\epsilon, x_1,
\ldots,x_n]\,$ and we let $J$ denote the intersection of $I$ with this
subring of $\cc(\epsilon)[x_1^{\pm}, \ldots,x_n^{\pm}]$.  Generators
of $J$ are computed from generators of $I$ by clearing denominators
and saturating with respect to $\epsilon$. The tropical variety of
$I'$ is related to the tropical variety of $J$ as follows.

\begin{lemma} \label{northern}
  A vector $w \in \rr^n$ lies in the polyhedral complex
  $\mathcal{T}(I')$ if and only if the vector $(1,w) \in \rr^{n+1}$
  lies in the polyhedral fan $\mathcal{T}(J)$.
\end{lemma}

Thus the tropical variety $\mathcal{T}(I')$ equals the restriction of
$\mathcal{T}(J)$ to the northern hemisphere of the $n$-sphere.  Note
that if $I$ is a prime ideal then so are $I'$ and $J$. Einsiedler,
Kapranov and Lind \cite{EKL} have shown that if $I'$ is prime, then
$\mathcal{T}(I')$ is connected.  Our connectivity results in Section 3
(which use the result of \cite{EKL}) imply the following
result which was conjectured in \cite{EKL}.

\begin{theorem}
  If $I$ is an ideal in $\,\cc \{ \! \{ \epsilon \} \! \}[x_1,
  \ldots,x_n]$ whose radical is prime of dimension $d$, then the
  tropical variety $\mathcal{T}(I)$ is a pure $d$-dimensional
  polyhedral complex which is connected in codimension one.
\end{theorem}

On the algorithmic side, we conclude that the polyhedral complex
$\mathcal{T}(I')$ can be computed by restricting the flip algorithm of
Section 4 to maximal cones in the fan $\mathcal{T}(J)$ which
intersect the open northern hemisphere in $\rr^{n+1}$.

\section{Algorithmic Problems and Tropical Bases}
\label{sec:algorithmic problems}

For all algorithms in this paper we fix the
ambient ring to be the polynomial ring over the complex numbers,
$\cc[\xx] := \cc[x_1, \ldots, x_n]$. The most basic
computational problem in tropical geometry is the following:

\begin{problem} \label{prob:tropical_prevariety}
  Given a finite list of polynomials $f_1, \ldots, f_r \in \cc[\xx]$,
  compute the tropical prevariety $\mathcal{T}(f_1) \cap \cdots \cap
  \mathcal{T}(f_r)$ in $\rr^n$~.
\end{problem}

The geometry of this problem is best understood by considering the
Newton polytopes ${\rm New}(f_1), \ldots,{\rm New}(f_r)$ of the given
polynomials. By definition, ${\rm New}(f_i)$ is the convex hull in
$\rr^n$ of the exponent vectors which appear with non-zero coefficient
in $f_i$. The tropical hypersurface $\mathcal{T}(f_i)$ is the
$(n-1)$-skeleton of the inner normal fan of the polytope ${\rm
  New}(f_i)$. Our problem is to intersect these normal fans. The
resulting tropical prevariety can be a fairly general polyhedral fan.
Its maximal cones may have different dimensions.

The tropical variety of an ideal $I$ in $\cc[\xx]$ is the set
$\trop(I) := \bigcap_{f \in I} \trop(f)$.  Equivalently, $\mathcal
T(I) = \{ w \in \mathbb R^n \,:\, \textup{in}_w(I) \,\textup{does
  not contain a monomial} \}$ where $\textup{in}_w(I) := \langle \,
\textup{in}_w(f) \, : \, f \in I \, \rangle$ is the initial ideal of
$I$ with respect to $w$. Bieri and Groves \cite{BG}
proved that $\mathcal{T}(I)$ is a $d$-dimensional fan when $d$ is the
Krull dimension of $\cc[\xx]/I$. The fan is pure if $I$ is
unmixed.  In Section 3 we shall prove that $\mathcal{T}(I)$ is
connected in codimension one if $I$ is prime. 

We first note that it suffices to devise algorithms for computing
tropical varieties of homogeneous ideals. Let $^hI \subset \cc[x_0,
x_1, \ldots, x_n]$ be the homogenization of an ideal $I$ in
$\cc[\xx]$ and $^hf$ the homogenization of $f \in \cc[\xx]$.

\begin{lemma}\label{lemma: homogeneous}
  Fix an ideal $I \subset \cc[\xx]$ and a vector $w \in \rr^n$.
  The initial ideal ${\rm in}_w(I)$ contains a monomial if and
  only if ${\rm in}_{(0,w)}(^hI)$ contains a monomial.
\end{lemma}

\begin{proof}
  Suppose $\xx^\u \in {\rm in}_w(I)$. Then $\xx^\u={\rm in}_w(f)$ for
  some $f \in I$. The $(0,w)$-weight of a term in $^hf$ equals the
  $w$-weight of the corresponding term in $f$. Hence ${\rm
    in}_{(0,w)}(^hf)=x_0^a\xx^\u \in {\rm in}_{(0,w)}(^hI)$ where $a$
  is some non-negative integer.
  
  Conversely, if $\xx^\u \in {\rm in}_{(0,w)}(^hI)$ then $\xx^\u =
  {\rm in}_{(0,w)}(f)$ for some $f \in {}^hI$.  Substituting $x_0=1$
  in $f$ gives a polynomial in $I$. The $(0,w)$-weight of any term in
  $f$ equals the $w$-weight of the corresponding term in
  $\,f|_{x_0=1}$. Since ${\rm in}_{(0,w)}(f)$ is a monomial, only one
  term in $f$ has minimal $(0,w)$-weight. This term cannot be canceled
  during the substitution. Hence it lies in ${\rm in}_w(I)$.
\end{proof}

Our main goal in this paper is to solve the following problem.

\begin{problem} \label{prob:tropical_variety}
  Given a finite list of homogeneous polynomials $f_1, \ldots, f_r \in
  \cc[\xx]$, compute the tropical variety $\,\mathcal{T}(I)$ of their
  ideal $\,I = \langle f_1,\ldots,f_r \rangle$.
\end{problem}

It is important to note that the two problems stated so far are of a
fundamentally different nature. Problem \ref{prob:tropical_prevariety}
is a problem of polyhedral geometry. It involves only polyhedral
computations: no algebraic computations are required. Problem
\ref{prob:tropical_variety}, on the other hand, combines the
polyhedral aspect with an algebraic one. To solve Problem
\ref{prob:tropical_variety} we must perform algebraic operations
(e.g.~Gr\"obner bases) with polynomials. In
Problem~\ref{prob:tropical_prevariety} we do not assume that the input
polynomials $f_1, \ldots, f_r$ are homogeneous as the polyhedral
computations can be performed easily without this assumption.

\begin{proposition} \label{HomogTwoWays}
Let $I$ be an ideal in $\cc[\xx]$ and let $w \in \rr^n$. The following are equivalent:
\begin{enumerate}
\item The ideal $I$ is $w$-homogeneous; \emph{i.e.} $I$ is generated by a 
set $S$ of $w$-homogeneous polynomials, meaning that $\myinit_w(f)=f$ for all $f\in S$. 
\item The initial ideal $\myinit_w(I)$ is equal to $I$.
\end{enumerate}
\end{proposition}

\begin{proof}
  If $I$ has a $w$-homogeneous generating set then
  $I \subseteq \myinit_w(I)$. Any maximal $w$-homogeneous component of
  $f\in I$ is in $I$. In particular $\myinit_w(f)\in I$.  Conversely,
  the ideal $\myinit_w(I)$ is generated by $w$-homogeneous elements by
  definition so, if $I=\myinit_w(I)$, then $I$ is generated by
  $w$-homogeneous elements.
\end{proof}

%

The set of $w \in \rr^n$ for which the above equivalent conditions
hold is a vector subspace of $\rr^n$. Its dimension is called the {\em
  homogeneity} of $I$ and is denoted ${\rm homog}(I)$. This space is
contained in every cone of the fan $\mathcal{T}(I)$ and can be
computed from the Newton polytopes of the polynomials that form any
reduced Gr\"obner basis of $I$. Passing to the quotient of $\rr^n$
modulo that subspace and then to a sphere around the origin,
$\mathcal{T}(I)$ can be represented as a polyhedral complex of
dimension $\, n - {\rm codim}(I) - {\rm homog}(I) - 1 = {\rm dim}(I) -
{\rm homog}(I) - 1$. Here ${\rm codim}(I)$ and ${\rm dim}(I)$ are the
codimension and dimension of $I$. In what follows, $\mathcal{T}(I)$ is
always presented in this way, and every ideal $I$ is presented by a
finite list of generators together with the three numbers $ n$, $ {\rm
  dim}(I) $ and $ {\rm homog}(I)$.

\begin{example}
\label{sym4x4}
Let $I$ denote the ideal which is generated by the
$3 \times 3$-minors of a symmetric $4 \times 4$-matrix
of unknowns. This ideal has $n = 10$,
${\rm dim}(I) = 7$ and ${\rm homog}(I) = 4 $.
Hence $\mathcal{T}(I)$ is a two-dimensional
polyhedral complex. We regard $\mathcal{T}(I)$ as the
tropicalization of the  secant variety of the
Veronese threefold in $\pp^9$, i.e., the variety of
symmetric $4 \times 4$-matrices of rank $\leq 2$,
Applying our  {\tt Gfan} implementation
(see Example \ref{3minors of nbynsymmetric}),
we find that $\mathcal{T}(I)$ is a simplicial complex consisting
of $75$ triangles, $75$ edges and $20$ vertices.
 \qed
\end{example}
 
Our next problem concerns tropical bases.  A finite set $\{f_1
,\ldots, f_t\} $ is a tropical basis of $I$ if $ \,\langle f_1,
\ldots, f_t \rangle = I \,\, \textup{and} \,\, \mathcal{T}(I) =
\mathcal{T}(f_1) \cap \cdots \cap \mathcal{T}(f_t).$
 
\begin{problem} \label{prob:find tropical basis}
Compute a tropical basis of a given ideal $\,I \subset \cc[\xx]$. 
\end{problem}

A priori, it is not clear that every ideal $I$ has a finite tropical basis,
but we shall prove this below. First, here is one case
where this is easy:

\begin{example} \label{rem:principal} If $\, I = \langle f \rangle \,$ is a
principal ideal, then $\{f\}$ is a tropical basis. \qed
\end{example}

In \cite{SS} it was claimed that any universal Gr\"obner basis of $I$ is
a tropical basis. Unfortunately, this claim is false as the
following example shows.

\begin{example} \label{ex:universalGB}
  Let $I$ be the intersection of the three linear ideals $\langle x+y, z
  \rangle$, $\langle x+z, y \rangle$, and $\langle y+z, x \rangle$
  in $\cc[x,y,z]$.
  Then $I$ contains the monomial $xyz$, so $\TT(I)$ is empty.  
  A minimal universal Gr\"obner basis of $I$ is
$$ \mathcal{U} \quad = \quad \{\,x+y+z,\,
 x^2y+xy^2 , \,y^2z+yz^2 ,\,x^2z+xz^2 \, \}, $$
and the intersection of the four corresponding tropical surfaces
in $ \rr^3$ is the line $w_1 = w_2 = w_3 $.
 Thus $\mathcal{U}$ is not a tropical basis of $I$.  \qed
\end{example}

We now prove that every ideal $I \subset \cc[\xx]$ has a tropical
basis.  By Lemma~\ref{lemma: homogeneous}, one tropical basis of a
non-homogeneous ideal $I$ is the dehomogenization of a tropical basis
for $^hI$. Hence we shall assume that $I$ is a homogeneous ideal.

Tropical bases can be constructed from the {\em Gr\"obner fan} of $I$
(see \cite{MR}, \cite{GBCP}) which is a complete finite rational
polyhedral fan in $\rr^n$ whose relatively open cones are in bijection
with the distinct initial ideals of $I$. Two weight
vectors $w, w' \in \rr^n$ lie in the same relatively open cone of the
Gr\"obner fan of $I$ if and only if ${\rm in}_w(I) = {\rm
  in}_{w'}(I)$. The closure of this cell, denoted by $C_w(I)$, is
called a {\em Gr\"obner cone} of $I$. The $n$-dimensional Gr\"obner cones
are in bijection with the reduced Gr\"obner bases, or equivalently,
the monomial initial ideals of $I$.  Every Gr\"obner cone of $I$ is a
face of at least one $n$-dimensional Gr\"obner cone of $I$. If ${\rm
  in}_w(I)$ is not a monomial ideal, then we can refine $w$ to
$\prec_w$ by breaking ties in the partial order induced by $w$ with a
fixed term order $\prec$ on $\cc[\xx]$. Let $\G_{\prec_w}(I)$ denote
the reduced Gr\"obner basis of $I$ with respect to $\prec_w$. The
Gr\"obner cone of $\G_{\prec_w}(I)$, denoted by $C_{\prec_w}(I)$, is
an $n$-dimensional Gr\"obner cone that has $C_w(I)$ as a
face.  The tropical variety $\trop(I)$ consists of all Gr\"obner cones
$C_w(I)$ such that ${\rm in}_w(I)$ does not contain a monomial. From the description of $\mathcal{T}(I)$ as $\,\bigcap_{f \in I} \mathcal{T}(f)\,$ it is clear that $\mathcal{T}(I)$ is closed. Thus we deduce that $\mathcal{T}(I)$ is a closed subfan of the Gr\"obner fan. This
endows the tropical variety $\TT(I)$ with the structure of a polyhedral fan.

\begin{theorem}  \label{thm:tropical bases}
  Every ideal $I \subset \cc[\xx]$ has a tropical basis.
\end{theorem}

\begin{proof}  Let $\mathcal{F}$ be any finite generating set of $I$ which is not a tropical basis. Pick a
  Gr\"obner cone $C_w(I)$ whose relative interior intersects
  $\,\cap_{f \in \mathcal{F}} \TT(f)\,$ non-trivially and whose
  initial ideal $\,{\rm in}_w(I)\,$ contains a monomial $\xx^{\bf m}$.
  Compute the reduced Gr\"obner basis $\G_{\prec_w}(I)$ for a
  refinement $\prec_w$ of $w$, and let $h$ be the normal form of
  $\xx^{\bf m}$ with respect to $\G_{\prec_w}(I)$.  Let $f := \xx^{\bf
    m} - h$. Since the normal form of $\xx^{\bf m}$ with respect to
  $\G_\prec({\rm in}_w(I))=\{{\rm in}_w(g):g\in\G_{\prec_w}(I)\}$ is
  $0$ and $h$ is the normal form of $\xx^{\bf m}$ with respect to
  $\G_{\prec_w}(I)$, every monomial occurring in $h$ has higher
  $w$-weight than $\xx^{\bf m}$.  Moreover, $h$ depends only on the
  reduced Gr\"obner basis $\G_{\prec_w}(I)$ and is independent of the
  particular choice of $w$ in $C_w(I)$. Hence for any $w'$ in the
  relative interior of $C_w(I)$, we have $\xx^{\bf m} = {\rm
    in}_{w'}(f)$. This implies that the polynomial $f:= \xx^{\bf m} -
  h$ is a {\em witness} for the cone $C_w(I)$ not being in the
  tropical variety $\trop(I)$.
  
  We now add the witness $f$ to the current basis $\mathcal{F}$ and
  repeat the process. Since the Gr\"obner fan has only finitely many
  cones, this process will terminate after finitely many steps. It
  removes all cones of the Gr\"obner fan which violate the condition
  for $\mathcal{F}$ to be a tropical basis.
\end{proof}

We next show that tropical bases can be very
large even for linear ideals.  Let $I $ be the ideal in $\cc[\xx]$
generated by $d$ linear forms $\sum_{j=1}^{n} a_{ij} x_j$
where $i =1, \ldots, d$ and $(a_{ij})$ is an integer $d \times n$
matrix of rank $d$.  The tropical variety $\TT(I)$ depends only on the
matroid associated with $I$, and it is known as the {\em Bergman fan}
of that matroid. The results on the Bergman fan proved in \cite{AK, cbms}
imply that the circuits in $I$ form a tropical basis.
A {\em circuit} of $I$ is a non-zero linear polynomial $f \in I$ of
minimal support.  The following result answers the
question which was  posed in \cite[\S 5]{SS1}.

\begin{theorem} For any $1 \leq d \leq n$, there is a linear ideal
  $I$ in $\cc[x_1,\ldots,x_n]$ such that any tropical basis of linear
  forms in $I$ has size at least $\frac{1}{n-d+1} {n \choose d}$.
\end{theorem}

\begin{proof} Suppose that all $d \times d$-minors of the coefficient
  matrix $(a_{ij})$ are non-zero. Equivalently, the matroid of $I$ is
  uniform.  There are ${n \choose n-d+1}$ circuits in $I$, each
  supported on a different $(n-d+1)$-subset of $\{x_1,\ldots,x_n\}$.
  Since the circuits form a tropical basis of $I$ and each circuit has
  support of size $n-d+1$, the tropical variety $\TT(I)$
  consists of all vectors $w \in \rr^n$ whose smallest $d+1$
  components are equal. The latter condition is necessary and
  sufficient to ensure that no single variable in a circuit becomes 
  the initial form of the circuit with respect to $w$.
Consider any vector $w \in \rr^n$ satisfying
$$\,w_{i_1} = w_{i_2} = \cdots = w_{i_d} \,< \, {\rm min}\bigl( w_{j}
\,: \, j \in \{ 1,\ldots,n\} \backslash \{i_1,i_2,\ldots,i_d\} \bigr)
.$$
Since $w \not\in \TT(I)$, any tropical basis of linear forms
in $I$ contains an $f$ such that $\,{\rm in}_w(f) \in
\{x_{i_1}, \ldots,x_{i_d}\}$. This implies that $f$ is
one of the $d$ circuits whose support contains the $n-d$ variables
$\,x_j \,$ with $ j \not\in \{i_1,\ldots,i_d\} $.  The support of each
circuit has size $n-d+1$, hence contains $n-d+1$ distinct
$(n-d)$-subsets. There are ${n \choose d}$ $\,(n-d)$-subsets of
$\{x_1,\ldots,x_n\}$ to be covered.  Hence any tropical basis
consisting of linear forms has size at least $\frac{1}{n-d+1}
{n \choose d}$.
\end{proof}

\begin{example}
Let $d=3, n=5$. The Bergman fan $\TT(I)$ corresponds to 
the line in tropical projective $4$-space which consists of
the five rays in the coordinate directions.
We have $\,\frac{1}{n-d+1} {n \choose d}\,=\,10/3$. Hence
this line is not a complete intersection
of three tropical hyperplanes, but it requires four. \qed
\end{example}

\section{Transversality and Connectivity}

In this section we assume that $I$ is a prime ideal of dimension $d$
in $\cc[x_1,\ldots,x_n]$. Then its tropical variety $\mathcal{T}(I)$
is called {\em irreducible}. It is a subfan of the Gr\"obner fan of
$I$ and, by the Bieri-Groves Theorem \cite{BG, cbms}, all facets of
$\mathcal{T}(I)$ are cones of dimension $d$.  A cone of dimension
$d-1$ in $\TT(I)$ is called a {\em ridge} of the tropical variety
$\mathcal{T}(I)$.  A {\em ridge path} is a sequence of facets $\,
F_1,F_2, \ldots, F_k \,$ such that $F_{i} \cap F_{i+1}$ is a ridge for
all $\, i \in \{1,2,\ldots,k-1\}$.  Our objective is to prove the
following result, which is crucial for the algorithms.

\begin{theorem} \label{codimone}
Any irreducible tropical variety $\mathcal{T}(I)$ is connected
in codimension one, i.e., any two facets are connected
by a ridge path.
\end{theorem}

The proof of this theorem will be based on the following
important lemma.

\begin{lemma} \label{intersect}
{\rm (Transverse Intersection Lemma)} \hfill \break
Let $I$ and $J$ be ideals in $\cc[x_1,\ldots,x_n]$
whose tropical varieties $\mathcal{T}(I)$
and $\mathcal{T}(J)$ meet transversally at
a point $w \in \rr^n$. Then $w \in \mathcal{T}(I+J)$.
\end{lemma}

By ``meet transversely'' we mean that
 if $F$ and $G$ are the cones of
$\mathcal{T}(I)$ and $\mathcal{T}(J)$ which contain $w$ 
in their relative interior,
 then $\rr F + \rr G=\rr^n$.

This lemma implies that any
transverse intersection of tropical varieties is a tropical variety.
In particular, any transverse intersection of tropical 
hypersurfaces is a tropical variety, and such a
tropical variety is defined by an ideal which is
a complete intersection in the commutative algebra sense.

\begin{corollary}
For any two ideals $I$ and $J$
 in $\cc[x_1,\ldots,x_n]$ we have
$$ \TT(I + J) \quad \subseteq \quad \TT(I) \,\cap \, \TT(J). $$
Equality holds if the latter intersection is transverse at every point
except the origin and the two fans meet in at least one point other
than the origin.
\end{corollary}

\begin{proof}
  We have $\TT(I) \cap \TT(J) = \bigcap_{f \in I} \TT(f) \cap
  \bigcap_{f \in J} \TT(f)=\bigcap_{f \in I \cup J} \TT(f)$. Clearly,
  this contains $\TT(I+J) = \bigcap_{f \in I+J} \TT(f)$. If $\TT(I)$
  and $\TT(J)$ intersect transversally and $w$ is a point of
  $\mathcal{T}(I) \cap \mathcal{T}(J)$ other than the origin then the
  preceeding lemma tells us that $w \in \mathcal{T}(I+J)$. Thus
  $\mathcal{T}(I+J)$ contains every point of
  $\mathcal{T}(I) \cap \mathcal{T}(J)$ except possibly the origin. In
  particular, $\mathcal{T}(I+J)$ is not empty. Every nonempty fan
  contains the origin, so we see that the origin is in
  $\mathcal{T}(I+J)$ as well.
\end{proof}

We first derive Theorem \ref{codimone} from Lemma \ref{intersect},
which will be proved later. We must at this point address an annoying
technical detail. The subset $\mathcal{T}(I) \subset \rr^n$ depends
only on the ideal $I \cc[\xx^{\pm 1}]$ generated by $I$ in the Laurent
polynomial ring $\cc[x_1^{\pm}, \ldots, x_n^{\pm 1}]$. (This is easy
to see: if $I_1$ and $I_2$ generate the same ideal in $\cc[\xx^{\pm 1}]$
and $w \not \in \mathcal{T}(I_1)$ then there is a polynomial $f \in
I_1$ such that $\mathrm{in}_w (f)$ is a monomial. There is some monomial
$m$ such that $m f \in I_2$, then $\mathrm{in}_w (mf)$ is a monomial
and $w \not \in \mathcal{T}(I_2)$.)  From a theoretical perspective
then, it would be better to directly work with ideals in $\cc[\xx^{\pm 1}]$.
One reason is the availability the symmetry group ${\mathrm GL}_n(\zz)$ of
the multiplicative group of monomials. The action of this group
transforms $\mathcal{T}(I)$ by the obvious action on $\rr^n$. This
symmetry will prove invaluable for simplifying the arguments in this
section. Therefore, in this section, we will work with ideals in
$\cc[\xx^{\pm 1}]$.     Computationally, however, it is much better to deal with
ideals in $\cc[\xx]$ as it is for such ideals that Gr\"obner basis
techniques have been developed and this is the approach we take in the
rest of the paper.

 Note that, if $I \subset \cc[\xx]$ is prime then
so is the ideal it generates in $\cc[\xx^{\pm 1}]$. We will signify an
application of the ${\mathrm GL}_n(\zz)$ symmetry by the phrase
``making a multiplicative change of variables''. The polyhedral
structure on $\mathcal{T}(I)$ induced by the Gr\"obner fan of $I$ may
change under a multiplicative change of variables of $I
\cc[\xx^{\pm}]$ in $\cc[\xx^{\pm 1}]$, but all of the properties of
$\mathcal{T}(I)$ that are of interest to us depend only on the
underlying point set.

\medskip

\noindent {\sl Proof of Theorem \ref{codimone}.}
As discussed, we replace $I$ by the ideal it generates in
$\cc[\xx^{\pm 1}]$ and, by abuse of notation, continue to denote this
ideal as $I$.  The proof is by induction on $d = {\rm
  dim}(\mathcal{T}(I))$.  If $d \leq 1$ then the statement is
trivially true. We now explain why the result holds for $d=2$.  By a
multiplicative change of coordinates, it suffices to check that
$\mathcal{T}(I) \cap \{ x_n = 1 \}$ is connected.
Let $K$ be the Puiseux series field over $\cc$. Let $I' \subset K[x_1,
\ldots, x_{n-1}]$ be the prime ideal generated by $I$ via the
inclusion $\cc[x_n] \to K$.  By Lemma \ref{northern}, the tropical
variety of $I'$ is $\mathcal{T}(I) \cap \{ x_n = 1 \}$.  In \cite{EKL}
it was shown that the tropical variety of $I'$ is connected whenever
$I'$ is prime.  We conclude that $\mathcal{T}(I) \cap \{ x_n = 1 \}$
is connected, so our result holds for $d=2$.

We now suppose that $d \geq 3$. Let $F$ and $F'$ be facets of
$\mathcal{T}(I)$.  We can find
$$
H \quad = \quad \bigl\{\,(u_1,\ldots,u_n) \in \rr^n \,:\, a_1 u_1 +
\cdots + a_n u_n = 0 \,\bigr\} $$
such that $a_1,\ldots,a_n$ are
relatively prime integers, both $H \cap F$ and $H \cap F'$ are cones
of dimension $d-1$, and $H$ intersects every cone of $\mathcal{T}(I)$
except for the origin transversally. To see this, select rays $w$
and $w'$ in the relative interiors of $F$ and $F'$. By perturbing $w$
and $w'$ slightly, we may arrange that the span of $w$ and $w'$ does
not meet any ray of $\mathcal{T}(I)$ -- here it is important that $d
\geq 3$. Now, taking $H$ to be the span of $w$, $w'$ and a generic $(n-3)$-plane, we get that $H$ also does not contain any ray of $\mathcal{T}(I)$ and hence does not contain any positive dimensional face of $\mathcal{T}(I)$. So $H$ is transverse to $\mathcal{T}(I)$ everywhere except at the origin.
 Since $H \cap F$ and $H \cap F'$ are
positive-dimensional (as $d \geq 2$), the hyperplane $H$ does intersect
$\mathcal{T}(I)$ at points other than just the origin. The hyperplane
$H$ is the tropical hypersurface of a binomial, namely, $H =
\mathcal{T}(\langle f_u \rangle)$, where
$$ f_u \quad = \quad \prod_{i: a_i > 0} (u_i x_i)^{a_i}
\, - \,  \prod_{j : a_j < 0} (u_j x_j)^{-a_j}, $$
and $\,u = (u_1,u_2,\ldots,u_n)\,$ is an arbitrary
point in the algebraic torus  $(\cc^*)^n$.
Our transversality assumption regarding $H$ and
Lemma \ref{intersect} imply that
\begin{equation}
\label{HcT}
 H \,\cap \, \mathcal{T}(I) \quad = \quad
\mathcal{T}(\langle f_u \rangle) \,\cap\,
\mathcal{T}(I) \quad =  \quad
\mathcal{T} \bigl( I + \langle f_u \rangle \bigr).
\end{equation}
Since $I$ is prime of dimension $d$,
and $f_u \not\in I$, the ideal $\,I + \langle f_u \rangle\,$
has dimension $d-1$ by Krull's Principal Ideal Theorem \cite[Theorem 10.1]{Eis}. If
$\,I + \langle f_u \rangle\,$ were a prime ideal
then we would be done by induction.
Indeed, this would imply that there is
a ridge path between the facets $H \cap F$ and
$H \cap F'$ in the
$(d-1)$-dimensional tropical variety
$(\ref{HcT})$. Since $d\geq 3$, the $(d-1)$- and $(d-2)$-dimensional faces of $H \cap \mathcal{T}(I)$ arise
uniquely
 from the intersections of $H$ with
$d$- and $(d-1)$-dimensional faces of $\mathcal{T}(I)$.
Hence this path is also a ridge path considered as a path in $\mathcal{T}(I)$.

Let $V(J)$ denote the subvariety of the algebraic torus
$  (\cc^*)^n$ defined by an ideal
$J \subset \cc[x_1^{\pm 1}, \ldots, x_n^{\pm 1}]$.
 The tropical variety in (\ref{HcT}) depends only on
 the  subvariety of $(\cc^*)^n$ defined by our ideal
$\, I + \langle f_u \rangle$. This subvariety is
 \begin{equation}
 \label{tropvar}
 V \bigl(I + \langle f_u \rangle \bigr) \quad = \quad
  V(I) \,\cap \, V(f_u) \quad  = \quad
V(I) \,\cap \, u^{-1} \cdot V(f_{\bf 1}) .
\end{equation}
Here ${\bf 1}$ denotes the identity element of $(\cc^*)^n$.
For generic choices of the group element
$u \in (\cc^*)^n$, the intersection (\ref{tropvar}) is
an irreducible subvariety of dimension $d-1$ in $(\cc^*)^n$. This follows
from Kleiman's version of Bertini's Theorem \cite[Theorem III.10.8]{Har},
applied to the algebraic group $(\cc^*)^n$. 
 Hence (\ref{HcT}) is indeed
an irreducible tropical variety of dimension $d-1$,
defined by the prime ideal $I  + \langle f_u \rangle$.
This completes the proof by induction.
\qed

\medskip

\noindent {\sl Proof of Lemma \ref{intersect}:}
Again, we replace $I \subset \cc[\xx]$ by the ideal it generates in
$\cc[\xx^{\pm 1}]$. Let $F$ be the cone of $\mathcal{T}(I)$ which contains
$w$ in its relative interior and  $G$
the cone of $\mathcal{T}(J)$ which contains
$w$ in its relative interior.
Our hypothesis is that $F$ and $G$ meet transversally at $w$, that is,
$$ \rr F \, + \, \rr G \quad = \quad \rr^n  . $$

We claim that the ideal ${\rm in}_w(I)$ is homogeneous with respect to
any weight vector $v\in\rr F$ or, equivalently (see Proposition \ref{HomogTwoWays}), that ${\rm in}_v({\rm
in}_w(I))={\rm in}_w(I)$. According to Proposition 1.13 in
\cite{GBCP}, for $\epsilon$ a sufficiently small positive number,
${\rm in}_{w+\epsilon v} (I) ={\rm in}_v({\rm in}_w (I))$. The vector
$w+\epsilon v$ is in the relative interior of $F$ so ${\rm
in}_{w+\epsilon v} (I)={\rm in}_w (I)$.  By the same argument, the
ideal ${\rm in}_w(J)$ is homogeneous with respect to any weight vector
in $\rr G$.


After a multiplicative change of variables
in $\,\cc[x_1^{\pm 1}, \ldots, x_n^{\pm 1}]\,$
we may assume that $w = e_1$,
$\,\rr \{e_1 ,e_2, \ldots,e_s\} \subseteq \rr F\,$ and
$\,\rr \{e_1,e_{s+1} , \ldots,e_n\} \subseteq \rr G$.
We change the notation for the variables as follows:
$$ t = x_1 , \,
y = (y_2,\ldots,y_s) = (x_2,\ldots,x_s),\,\,
z = (z_{s+1},\ldots,z_n) = (x_{s+1},\ldots,x_n) . $$
The homogeneity properties of the two initial ideals ensure that we
can pick generators $f_1(z), \ldots,$ $f_a(z)\,$ for  ${\rm in}_w(I)$
and generators $g_1(y), \ldots,g_b(y)$ for  ${\rm in}_w(J)$.
Since ${\rm in}_w(I)$ is not the unit ideal, the
Laurent polynomials $f_i(z)$ have a common
zero $Z = (Z_{s+1}, \ldots,Z_n) \in (\cc^*)^{n-s}$,
and likewise the Laurent polynomials
$g_j(y)$ have a common zero
$Y = (Y_2,\ldots,Y_s) \in (\cc^*)^{s-1}$.

Next we  consider the following general chain of inclusions
of ideals:
\begin{equation}
\label{cdotcap}
{\rm in}_w(I) \,\cdot {\rm in}_w(J) \,\subseteq \,
{\rm in}_w ( I \cdot J) \,\subseteq \,
{\rm in}_w( I \,\cap \,J) \, \subseteq\,
{\rm in}_w(I) \,\cap\, {\rm in}_w(J).
\end{equation}
The product of two ideals which are generated by
(Laurent) polynomials in disjoint sets of variables
equals the intersection of the two ideals.
Since the set of $y$-variables is disjoint from the
set of $z$-variables, it follows that the first ideal
in  (\ref{cdotcap}) equals the last ideal
in (\ref{cdotcap}). In particular, we conclude that
\begin{equation}
\label{goodnews}
{\rm in}_w( I \,\cap \,J) \quad= \quad
{\rm in}_w(I) \,\cap\, {\rm in}_w(J).
\end{equation}
We next claim that
\begin{equation}
\label{betternews}
{\rm in}_w( I \,+ \,J) 
\quad= \quad
{\rm in}_w(I) \,+ \, {\rm in}_w(J).
\end{equation}
The left  hand side is an ideal
which contains both ${\rm in}_w(I)$
and ${\rm in}_w(J)$, so it
contains their sum.
We must prove  that the right hand side
contains the left hand side.
Consider any element $f+g \in I+J$ where $f \in I$ and $g \in J$.
Let $\, f \,=\, f_0(y,z) + t \cdot f_1(t,y,z) \,$ and 
$\, g \,=\, g_0(y,z) + t \cdot g_1 (t,y,z)$.
We have the following representation for some
integer $a \geq 0$ and non-zero polynomial $h_0$:
$$ f+g \quad = \quad t^a \cdot h_0(y,z) + t^{a+1} \cdot h_1(t,y,z). $$
If $a = 0$ then we conclude
$$ {\rm in}_w(f+g) \,\, = \,\, h_0(y,z) \,\, = \,\,
 f_0(y,z) + g_0(y,z) \,\, \in \,\, {\rm in}_w(I)  + {\rm in}_w(J). $$
If $a \geq 1 $ then $\, f_0 \, = \, - g_0 \,$ lies in
$\,{\rm in}_w(I) \,\cap\, {\rm in}_w(J)$. In view of (\ref{goodnews}),
there exists $\, p \in \,I \cap J \,\,$ with
  $\, f_0 \, = \, - g_0 \, = \, {\rm in}_w(p)$.
Then $\, f+g = (f-p)+(g+p) \,$ and replacing
$f$ by $(f-p)/t$ and $g $ by $(g+p)/t$ puts us in the same
situation as before, but with $a$ reduced by $1$.
By induction on $a$, we conclude that
$\,{\rm in}_w(f+g)$ is in $\, {\rm in}_w(I)  + {\rm in}_w(J)$,
and  the claim (\ref{betternews}) follows.

For any constant $\,T \in \cc^*$, the vector
$\,(T,Y_2,\ldots,Y_s,Z_{s+1},\ldots,Z_n)\,$ is a common zero in
$(\cc^*)^n$ of the ideal (\ref{betternews}).  We conclude that $\,{\rm
  in}_w(I+J)\,$ is not the unit ideal, so it contains no monomial, and
hence $w \in \mathcal{T}(I+J)$.  \qed


\section{Algorithms} \label{sec:algorithms}

In this section we describe algorithms for solving the
computational problems raised in Section \ref{sec:algorithmic
problems}. The emphasis is on algorithms leading to a solution
of Problem \ref{prob:tropical_variety} for prime ideals, taking
advantage of Theorem \ref{codimone}.  Recall
that we only need to consider the case of homogeneous ideals in $\cc[\xx]$.

In order to state our algorithms we must first explain how polyhedral
cones and polyhedral fans are represented.  A polyhedral cone is
represented by a canonical minimal set of inequalities and equations.
Given arbitrary defining linear inequalities and equations, the task
of bringing these to a canonical form involves linear programming.
Representing a polyhedral fan requires a little thought. We are rarely
interested in all faces of all cones.

\begin{definition}
  A set $S$ of polyhedral cones in $\rr^n$ is said to {\em represent} a fan
  $\F$ in $\rr^n$ if the set of all faces of cones in $S$ is exactly
  $\F$.
\end{definition}

A representation may contain non-maximal cones, but each
cone is represented minimally by its canonical form.
A Gr\"obner cone $C_\andersomega(I)$ is represented by the pair
$(\G_{\prec_\andersomega}(\myinit_\andersomega(I)),\G_{\prec_\andersomega}(I))$
of marked reduced Gr\"obner bases, where $\prec$ is some globally
fixed term order. In a \emph{marked} Gr\"obner basis the initial terms
are distinguished. The advantage of using marked Gr\"obner bases is
that the weight vector $\andersomega$ need not
be stored -- we can deduce defining inequalities for its cone from the
marked reduced Gr\"obner bases themselves, see Example \ref{ex:hankel4}.
This is done as follows; see \cite[proof of Proposition 2.3]{GBCP}:

\begin{lemma}\label{lem: cone inequalities}
Let $I\subset \cc[\xx]$ be a homogeneous ideal, $\prec$ a term order
and $\andersomega\in\rr^n$ a vector. For any other vector
$\andersomega'\in\rr^n$:
$$\andersomega'\in C_{\andersomega}(I)  \quad
\Longleftrightarrow \quad \forall \, f\in\G_{\prec_\andersomega}(I): \, \myinit_{\andersomega}(\myinit_{\andersomega'}(f))=\myinit_{\andersomega}(f).$$
\end{lemma}

Our first two algorithms perform polyhedral computations, and they
solve Problem \ref{prob:tropical_prevariety}. By the {\em support} of
a fan we mean the union of its cones. Recall that, for a polynomial $f$,
the tropical hypersurface $\trop(f)$ is the union of
the normal cones of the edges of the Newton polytope
$\textup{New}(f)$.  The first algorithm computes these cones.

\begin{algorithm}Tropical Hypersurface\label{alg: hypersurface}\\
  {\bf Input:} $f\in \cc[\xx]$.\\
  {\bf Output:} A representation $S$ of a polyhedral fan whose support
  is $\trop(f)$.\\  
  $\{$\\
\TAB $S:=\emptyset$;\\
\TAB For every vertex $v\in\textup{New}(f)$\\
\TAB $\{$\\
\TAB \TAB Compute the normal cone $C$ of $v$ in $\textup{New}(f)$;\\
\TAB \TAB $S:=S\cup\{\textup{the facets of }C\}$;\\
\TAB  $\}$\\
  $\}$
\end{algorithm}

Let $\F_1$ and $\F_2$ be polyhedral fans in $\rr^n$. Their {\em common
refinement} is 
$$\F_1\wedge \F_2:=\{C_1\cap C_2\}_{(C_1,C_2)\in \F_1\times \F_2}.$$
To compute a common refinement we simply run through all pairs of
cones in the fan representations and bring their intersection to
canonical form. The canonical form makes it easy to remove duplicates.
\begin{algorithm}Common Refinement\label{alg: common refinement}\\
  {\bf Input:} Representations $S_1$ and $S_2$ for polyhedral fans $\F_1$ and $\F_2$.\\
  {\bf Output:} A representation $S$ for the common refinement  $\F_1\wedge\F_2$.\\
    $\{$\\
\TAB $S:=\emptyset$;\\
\TAB For every pair $(C_1,C_2)\in S_1\times S_2$\\
\TAB \TAB $S:=S\cup\{C_1\cap C_2\}$;\\
  $\}$
\end{algorithm}

If refinements of more than two fans are needed, Algorithm \ref{alg:
  common refinement} can be applied successively. Note that the
intersection of the support of two fans is the support of the fans'
common refinement. Hence Algorithm \ref{alg: common refinement} can be
used for computing intersections of tropical hypersurfaces. This solves Problem \ref{prob:tropical_prevariety}, but the output may be a highly redundant representation.

Recall (from the proof of Theorem \ref{thm:tropical bases}) that a
witness $f\in I$ is a polynomial which certifies $\,\trop(f)\cap
\textup{rel int}(C_\andersomega(I))=\emptyset$.  Computing witnesses
is essential for solving Problems \ref{prob:tropical_variety} and
 \ref{prob:find tropical basis}.  The first step of
constructing a witness is to check if the ideal
$\myinit_\andersomega(I)$ contains monomials, and, if so, compute one
such monomial.  The check for monomial containment can be implemented
by saturating the ideal with respect to the product of the variables
(cf.~\cite[Lemma 12.1]{GBCP}). Knowing that the ideal contains a
monomial, a simple way of finding one is to repeatedly reduce powers
of the product of the variables by applying the division algorithm
until the remainder is $0$.

\begin{algorithm}Monomial in Ideal\label{alg: monomial}\\
  {\bf Input:} A set of generators for an ideal $I\subset \cc[\xx]$. \\ 
  {\bf Output:} A monomial $m\in I$ if one exists, {\bf no} otherwise. \\
$\{$\\
\TAB If $((I:x_1\cdots x_n^\infty)\not=\langle 1\rangle)$ return {\bf no};\\
\TAB $m:=x_1\cdots x_n$;\\
\TAB While $(m\not\in I)$ $\,m:=m \cdot x_1\cdots x_n$;\\
\TAB Return $m$;\\
$\}$
\end{algorithm}

\begin{remark}
To pick the smallest monomial in $I$ with respect to a term order,
we first compute the largest monomial ideal contained in $I$ using \cite[Algorithm 4.2.2]{SST} and then pick the smallest monomial generator of this ideal.
\end{remark}

Constructing a witness from a monomial was already explained in the  proof
of Theorem \ref{thm:tropical bases}. We only state the input and
output of this algorithm.

\begin{algorithm}Witness\label{alg: witness}\\
  {\bf Input:} A set of generators for an ideal $I\subset \cc[\xx]$ and a vector 
  $\andersomega\in\rr^n$ with $\myinit_\andersomega(I)$ containing a monomial.\\
  {\bf Output:} A polynomial $f\in I$ such that the tropical hypersurface $\trop(
  f)$ and the relative interior of $C_\andersomega(I)$ have empty
  intersection.
\end{algorithm}

Combining Algorithm \ref{alg: monomial} and Algorithm \ref{alg:
  witness} with known methods (e.g.~\cite[Algorithm 3.6]{GBCP}) for
computing Gr\"obner fans, we can now compute the tropical variety
$\trop(I)$ and a tropical basis of $I$. This solves Problem
\ref{prob:tropical_variety} and Problem \ref{prob:find tropical
  basis}.  This approach is not at all practical, as shown in Section
6.

We will present a practical algorithm for computing $\trop(I)$ when
$I$ is prime. An ideal $I \subset \cc[\xx] $ is said to define a
\emph{tropical curve} if $\idim(I)=1+\ihom(I)$. Our problems are
easier in this case because a tropical curve consists of only finitely
many rays and the origin modulo the homogeneity space.

\begin{algorithm}Tropical Basis of a Curve \label{alg: tropical basis of curve}\\
{\bf Input:} A set of generators $\G$ for an ideal $I$ defining a tropical curve.\\
{\bf Output:} A tropical basis $\G'$ of $I$.\\
$\{$\\
\TAB Compute a representation $S$ of $\bigwedge_{g\in \G} \trop(g)$;\\
\TAB For every $C\in S$\\
\TAB $\{$\\
\TAB \TAB Let $\andersomega$ be a generic relative interior point in $C$;\\
\TAB \TAB If ($\myinit_\andersomega(I)$ contains a monomial)\\
\TAB \TAB \TAB then add a witness to $\G$ and restart the algorithm;\\
\TAB $\}$\\
\TAB $\G':=\G$;\\
$\}$
\end{algorithm}

\noindent
{\sl Proof of correctness. }  The algorithm terminates because $I$ has
only finitely many initial ideals and at least one is excluded in
every iteration.  If a vector $\andersomega$ passes the monomial test
(which verifies $w \in \trop(I)$) then $C$ has dimension $0$ or $1$
modulo the homogeneity space since we are looking at a curve and
$\andersomega$ is generic in $C$. Any other relative interior point of
$C$ would also have passed the monomial test.  (This property fails in
higher dimensions, when $\trop(I)$ is no longer a tropical curve).
Hence, when we terminate only points in the tropical variety are
covered by $S$. Thus $\G'$ is a tropical basis. \qed
  
  \smallskip
  
  In the curve case, combining Algorithms \ref{alg: hypersurface} and
  \ref{alg: common refinement} with Algorithm \ref{alg: tropical basis
    of curve} we get a reasonable method for solving Problem
  \ref{prob:tropical_variety}. This method is used as a subroutine in
  Algorithm \ref{alg: neighbours} below. In the remainder of this
  section we concentrate on providing a better algorithm for Problem
  \ref{prob:tropical_variety} in the case of a prime ideal. The idea
  is to use connectivity to traverse the tropical variety.

The next algorithm is an important subroutine for us.
We only specify the input and output. This algorithm is
one step in the Gr\"obner walk  \cite{collart}.

\begin{algorithm}Lift \label{alg: lift}\\
{\bf Input:} Marked reduced Gr\"obner bases $\G_{\prec'}(I)$ and 
$\G_{\prec_\andersomega}(\myinit_\andersomega(I))$ where
$\andersomega\in C_{\prec'}(I)$ is an unspecified vector and $\prec$ and
$\prec'$ are unspecified term orders.\\
{\bf Output:} The marked reduced Gr\"obner basis $\G_{\prec_\andersomega}(I)$.\\
\end{algorithm}

We now suppose that $I$ is a monomial-free prime ideal with
$d=\idim(I)$, and $\prec$ is a globally fixed term order.
 We first describe the local computations needed for a traversal of
the $d$-dimensional Gr\"obner cones contained in $\TT(I)$. 

\begin{algorithm} Neighbors \label{alg: neighbours}\\
  {\bf Input:} A pair $(\G_{\prec_\andersomega}(\myinit_\andersomega(I)),\G_{\prec_\andersomega}(I))$ such that $\myinit_\andersomega(I)$ is monomial-free and $C_\andersomega(I)$ has dimension $d$.\\
  {\bf Output:} The collection $N$ of pairs of the form $(\G_{\prec_{\andersomega'}}(\myinit_{\andersomega'}(I)),\G_{\prec_{\andersomega'}}(I))$ where one $\andersomega'$ is taken from the relative interior of each $d$-dimensional Gr\"obner cone contained in $\trop(I)$ that has a facet in common with $C_\andersomega(I)$.\\
\noindent
$\{$\\
\TAB $N:=\emptyset$;\\
\TAB Compute the set $\F$ of facets of $C_\andersomega(I)$;\\
\TAB For each facet $F\in \F$\\
\TAB $\{$\\
\TAB \TAB Compute the initial ideal $J:=\myinit_\u(I)$\\
\TAB \TAB \TAB where $\u$ is a relative interior point in $F$;\\
\TAB \TAB Use Algorithm \ref{alg: tropical basis of curve} and Algorithm \ref{alg: common refinement} to produce a relative\\
\TAB \TAB \TAB interior point $\v$ of each ray in the curve $\trop(J)$;\\
\TAB \TAB For each such $\v$\\
\TAB \TAB $\{$\\
\TAB \TAB \TAB Compute
$(\G_{\prec_{\v}}(\myinit_{\v}(J)),\G_{\prec_{\v}}(J))=(\G_{\prec_{\v_\u}}(\myinit_{\v}(J)),\G_{\prec_{\v_\u}}(J))$;\\
\TAB \TAB \TAB Apply Algorithm \ref{alg: lift} to $\G_{\prec_\andersomega}(I)$ and $\G_{\prec_{\v_\u}}(J)$ to get $\G_{\prec_{\v_\u}}(I)$;\\
\TAB \TAB \TAB $N:=N\cup\{(\G_{\prec_{\v_\u}}(\myinit_{\v}(J)),\G_{\prec_{\v_\u}}(I))\}$;\\
\TAB \TAB $\}$\\
\TAB $\}$\\
$\}$
\end{algorithm}

\begin{figure}
\epsfig{file=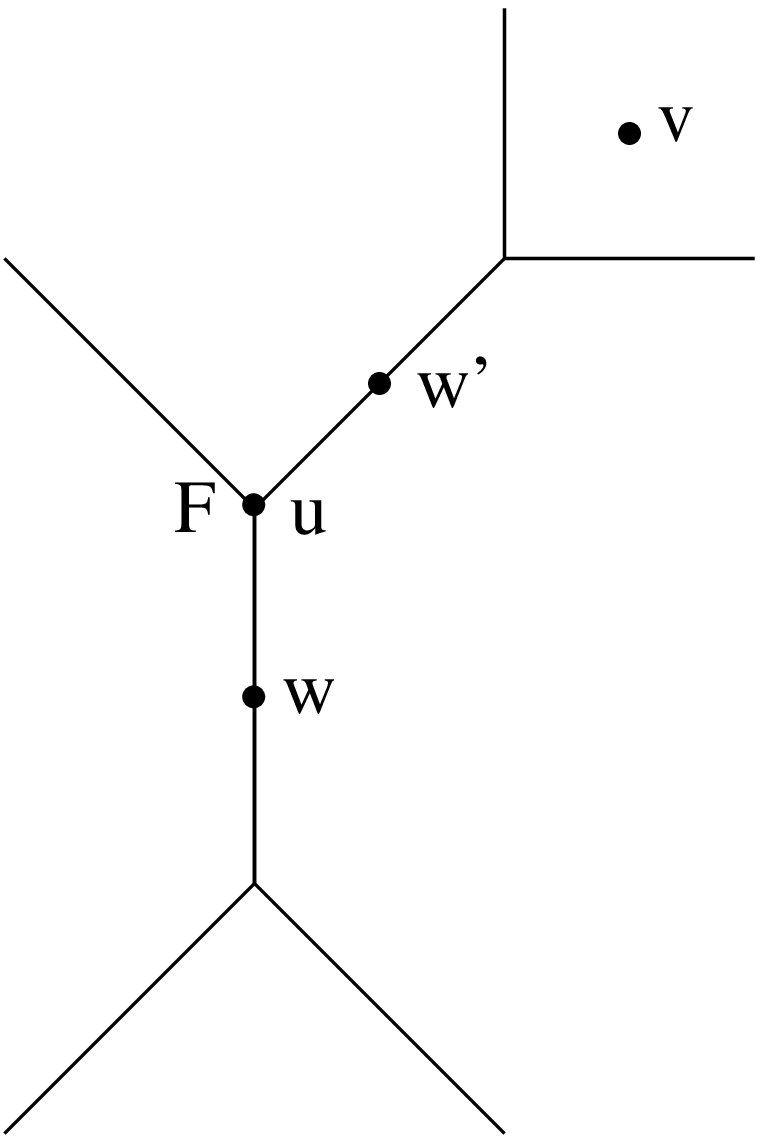,width=4cm}
\hspace{1cm}
\epsfig{file=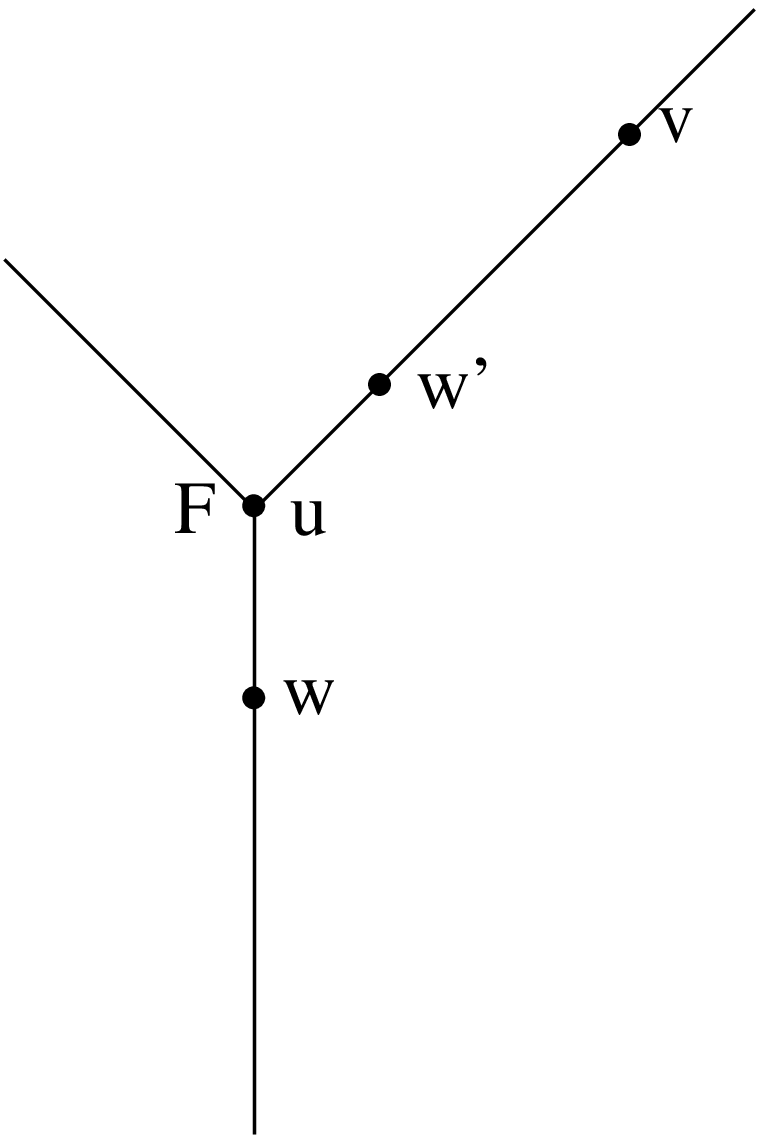,width=4cm}
\caption{A projective drawing of the situation in Algorithm \ref{alg: neighbours},
with $\trop(I)$  on the left and $\trop(\myinit_\u(I))$ on the right.}
\label{fig: neighbours}
\end{figure}

\noindent
{\sl Proof of correctness. }  
  Facets and relative interior points are computed
  using linear programming. Figure \ref{fig: neighbours} illustrates
  the choices of vectors in the algorithm.  The initial ideal
  $\myinit_\u(I)$ is homogeneous with respect to the 
  span of $F$. Hence its homogeneity space has dimension $d-1$. The
  Krull dimension of $\cc[\xx]/\myinit_\u(I)$ is $d$. Hence $\myinit_\u(I)$ defines a
  curve and $\trop(\myinit_\u(I))$ can be computed using
  Algorithm   \ref{alg: tropical basis of curve}. The identity
$\myinit_\v(\myinit_\u(I))=\myinit_{\u+\varepsilon\v}(I)$ for small
$\varepsilon>0$, see \cite[Proposition 1.13]{GBCP}, implies that we run through all the desired
$\myinit_{\andersomega'}(I)$ where $\andersomega'=\u+\varepsilon \v$
for small $\varepsilon>0$. The lifting step can be carried out since
$\u\in C_{\prec_\andersomega}(I)$.
 \qed
  
  \smallskip

\begin{algorithm}Traversal of an Irreducible Tropical Variety \label{alg: traverse}\\
  {\bf Input:} A pair $(\G_{\prec_\andersomega}(\myinit_\andersomega(I)),\G_{\prec_\andersomega}(I))$ such that $\myinit_\andersomega(I)$ is monomial free and $C_\andersomega(I)$ has dimension $d$.\\
  {\bf Output:}
  The collection $T$ of pairs of the form $(\G_{\prec_{\andersomega'}}(\myinit_{\andersomega'}(I)),\G_{\prec_{\andersomega'}}(I))$ where one $\andersomega'$ is taken from the relative interior of each $d$-dimensional Gr\"obner cone contained in $\trop(I)$. The union of all the $C_{\andersomega'}(I)$ is $\trop(I)$.\\
$\{$\\
\TAB $T:=\{(\G_{\prec_\andersomega}(\myinit_\andersomega(I)),\G_{\prec_\andersomega}(I))\}$;\\
\TAB $Old:=\emptyset$;\\
\TAB While ($T\not=Old$)\\
\TAB $\{$\\
\TAB \TAB $Old:=T$;\\
\TAB \TAB $T:=T\cup \textup{Neighbors}(T)$;\\
\TAB $\}$\\
$\}$
\end{algorithm}
\noindent
{\sl Proof of correctness. }  
By $\textup{Neighbors}(T)$ we mean the union of all the output of Algorithm
\ref{alg: neighbours} applied to all pairs in $T$.
The algorithm computes the connected component of the starting pair.
Since $I$ is a prime ideal, Theorem \ref{codimone} implies that the union of all the computed $C_{\andersomega'}(I)$ is $\trop(I)$.
 \qed
  
  \smallskip

To use Algorithm \ref{alg: traverse} we must know a starting
$d$-dimensional Gr\"obner cone contained in the tropical variety. One
inefficient method for finding one would be to compute the entire
Gr\"obner fan. Instead we currently use heuristics, which  are based
on the following probabilistic recursive algorithm:

\begin{algorithm} Starting Cone \label{alg: starting cone}\\
  {\bf Input:} A marked reduced Gr\"obner basis $\G$ for an ideal $I$ whose tropical variety is pure of dimension $d=\idim(I)$. A term order $\prec$ for tie-breaking.
  \\
  {\bf Output:} Two marked reduced Gr\"obner bases:
\begin{itemize}
\item One for an initial ideal $\myinit_{\andersomega'}(I)$ without monomials, where
the homogeneity space of $\myinit_{\andersomega'}(I)$ has dimension $d$. The term order is $\prec_{\andersomega'}$.
\item A marked reduced Gr\"obner basis for $I$ with respect to
  $\prec_{\andersomega'}$.
\end{itemize}
$\{$\\
\TAB If ($\idim(I) = \ihom(I)$)\\
\TAB \TAB Return $(\G_\prec(I),\G_\prec(I))$;\\
\TAB If not\\
\TAB $\{$\\
\TAB \TAB Repeat\\
\TAB \TAB $\{$\\
\TAB \TAB \TAB Compute a random reduced Gr\"obner basis of $I$;\\
\TAB \TAB \TAB Compute a random extreme ray $\andersomega$ of its Gr\"obner cone;\\
\TAB \TAB $\}$\\
\TAB \TAB Until ($\myinit_\andersomega(I)$ is monomial free);\\
\TAB \TAB Compute $\G_{\prec_\andersomega}(I)$;\\
\TAB \TAB $(\G_{Init},\G_{Full})$:= Starting Cone($\G_{\prec_\andersomega}(\myinit_\andersomega(I))$);\\
\TAB \TAB Apply Algorithm \ref{alg: lift} to $\G_{\prec_\andersomega}(I)$ and $\G_{Full}$\\
\TAB \TAB \TAB to get a marked reduced Gr\"obner basis $\G'$ for $I$;\\
\TAB \TAB Return $(\G_{Init},\G')$;\\
\TAB $\}$\\
$\}$
\end{algorithm}

\section{Software and Examples} \label{sec:examples}

We implemented the algorithms of Section~\ref{sec:algorithms} in the
software package \texttt{Gfan}~\cite{Gfan}.  \texttt{Gfan} uses the library
\texttt{cddlib}~\cite{cddlib} for polyhedral computations such as
finding facets and extreme rays of cones and bringing cones to
canonical form. We require our ideals to be generated by polynomials in $\qq[\xx]$. Exact arithmetic is done with the library
\texttt{gmp}~\cite{gmp}.  This is needed both for polyhedral
computations and for efficient
arithmetic in $\qq[\xx]$. In this section we illustrate the use of
\texttt{Gfan} in computing various tropical varieties.
 
\begin{example} \label{ex:hankel4}
We consider the prime ideal $\, I \subset \cc[a,b,c,d,e,f,g] \,$
which is generated by
the $3 \times 3$ minors of the generic  
\emph{Hankel matrix} of size $4 \times 4$:
$$ \left( \begin{array}{cccc}
a & b & c & d \\
b & c & d & e \\
c & d & e & f \\
d & e & f & g \\
\end{array} \right).$$
Its tropical variety  is a $4$-dimensional fan in $\rr^7$
with   $2$-dimensional homogeneity space.  Its combinatorics 
is given by the graph in Figure~\ref{fig: hankel}.   
To compute $\trop(I)$ in \texttt{Gfan}, we write the ideal
generators on a  file {\tt hankel.in}:

\begin{tiny}
\begin{verbatim}
% more hankel.in
{-c^3+2*b*c*d-a*d^2-b^2*e+a*c*e,-c^2*d+b*d^2+b*c*e-a*d*e-b^2*f+a*c*f,
-c*d^2+c^2*e+b*d*e-a*e^2-b*c*f+a*d*f,-d^3+2*c*d*e-b*e^2-c^2*f+b*d*f,
-c^2*d+b*d^2+b*c*e-a*d*e-b^2*f+a*c*f,-c*d^2+2*b*d*e-a*e^2-b^2*g+a*c*g,
-d^3+c*d*e+b*d*f-a*e*f-b*c*g+a*d*g,-d^2*e+c*e^2+c*d*f-b*e*f-c^2*g+b*d*g,
-c*d^2+c^2*e+b*d*e-a*e^2-b*c*f+a*d*f,-d^3+c*d*e+b*d*f-a*e*f-b*c*g+a*d*g,
-d^2*e+2*c*d*f-a*f^2-c^2*g+a*e*g,-d*e^2+d^2*f+c*e*f-b*f^2-c*d*g+b*e*g,
-d^3+2*c*d*e-b*e^2-c^2*f+b*d*f,-d^2*e+c*e^2+c*d*f-b*e*f-c^2*g+b*d*g,
-d*e^2+d^2*f+c*e*f-b*f^2-c*d*g+b*e*g,-e^3+2*d*e*f-c*f^2-d^2*g+c*e*g}
\end{verbatim}
\end{tiny}

We then run the command
\begin{verbatim}
gfan_tropicalstartingcone < hankel.in > hankel.start
\end{verbatim} 
which applies Algorithm~\ref{alg: starting cone} to produce a pair of marked
Gr\"obner bases. This represents
a maximal cone in $\trop(I)$,
as explained prior to Lemma
\ref{lem: cone inequalities}.
\begin{tiny}
\begin{verbatim}
% more hankel.start
{
c*f^2-c*e*g,
b*f^2-b*e*g,
b*e*f+c^2*g,
b*e^2+c^2*f,
b^2*g-a*c*g,
b^2*f-a*c*f,
b^2*e-a*c*e,
a*f^2-a*e*g,
a*e*f+b*c*g,
a*e^2+b*c*f}
{
c*f^2+e^3-2d*e*f+d^2*g-c*e*g,
b*f^2+d*e^2-d^2*f-c*e*f+c*d*g-b*e*g,
b*e*f+d^2*e-c*e^2-c*d*f+c^2*g-b*d*g,
b*e^2+d^3-2c*d*e+c^2*f-b*d*f,
b^2*g+c^2*e-b*d*e-b*c*f+a*d*f-a*c*g,
b^2*f+c^2*d-b*d^2-b*c*e+a*d*e-a*c*f,
b^2*e+c^3-2b*c*d+a*d^2-a*c*e,
a*f^2+d^2*e-2c*d*f+c^2*g-a*e*g,
a*e*f+d^3-c*d*e-b*d*f+b*c*g-a*d*g,
a*e^2+c*d^2-c^2*e-b*d*e+b*c*f-a*d*f}
\end{verbatim}
\end{tiny}
Using Lemma \ref{lem: cone inequalities} we can easily read off the canonical equations and equalities for the corresponding Gr\"obner cone $C_\andersomega(I)$.
 For example, the polynomials $cf^2-ceg$ and $cf^2+e^3-2def+d^2g-ceg$ 
 represent the equation
$$ \andersomega_c+2\andersomega_f \,\,= \,\, \andersomega_c+\andersomega_e+\andersomega_g$$
and the inequalities
$$
\andersomega_c+2\andersomega_f \,\,\leq \,\, {\rm min} \bigl\{
3\andersomega_e, \andersomega_d+\andersomega_e+\andersomega_f,
2\andersomega_d+\andersomega_g, \andersomega_c + \andersomega_e +
\andersomega_g \bigr\}. $$
At this point, we could run
Algorithm~\ref{alg: traverse} using the following command:
\begin{verbatim}
gfan_tropicaltraverse < hankel.start > hankel.out
\end{verbatim}

However, we can save computing time and get a better idea of the
structure of $\trop(I)$ by instructing \texttt{Gfan} to take advantage of
symmetries of $I$ as it produces cones.  The only symmetries that can
be used in \texttt{Gfan} are those that simply permute variables.  The
output will show which cones of $\trop(I)$ lie in the same orbit under
the action of the symmetry group we provide.

Our ideal $I$ is invariant under reflecting the
$4 \times 4$-matrix along the anti-diagonal. This 
reverses the variables $a,b,\ldots,g$.
To specify this permutation, we add
the following line to the bottom of the file
\texttt{hankel.start}:
\begin{verbatim}
{(6,5,4,3,2,1,0)}
\end{verbatim}
We can add more symmetries by listing them one after
another, separated by commas, inside the curly braces. \texttt{Gfan}
will compute and use the group generated by the set of permutations we
provide, and it will return an error if we input any permutation which
does not keep the ideal invariant. 
 
After adding the symmetries, we run the command
\begin{verbatim}
gfan_tropicaltraverse --symmetry < hankel.start > hankel.out
\end{verbatim}
to compute the tropical variety.  We show the output with some
annotations:

\vspace{.1in}
\begin{minipage}{2 in}
\begin{tiny}    
\begin{verbatim}
% more hankel.out

Ambient dimension: 7
Dimension of homogeneity space: 2 
Dimension of tropical variety: 4
Simplicial: true
Order of input symmetry group: 2
F-vector: (16,28)
\end{verbatim}
\end{tiny}
\end{minipage}
\begin{minipage}{2.7 in}
A short list of basic data: the dimensions of
the ambient space, of $\trop(I)$, and of its homogeneity space, and also 
the face numbers ($f$-vector) of $\trop(I)$ and the order of symmetry group specified in the input.   
\end{minipage}

\vspace{.1cm}

\begin{minipage}{2 in}
\begin{tiny}    
\begin{verbatim}
Modulo the homogeneity space:
{(6,5,4,3,2,-1,0),
 (5,4,3,2,1,0,-1)}
\end{verbatim}
\end{tiny}
\end{minipage}
\begin{minipage}{2.7 in}
A basis for the homogeneity space.  The rays are
considered in the quotient of $\rr^7$ modulo this $2$-dimensional subspace.
\end{minipage}

\begin{minipage}{2 in}
\begin{tiny}    
\begin{verbatim}
Rays: 
{0: (-1,0,0,0,0,0,0),
 1: (-5,-4,-3,-2,-1,0,0), 
 2: (1,0,0,0,0,0,0),
 3: (5,4,3,2,1,0,0),
 4: (2,1,0,0,0,0,0),
 5: (4,3,2,1,0,0,0),
 6: (0,-1,0,0,0,0,0),
 7: (6,5,4,3,2,0,0),
 8: (3,2,1,0,0,0,0),
 9: (0,0,-1,0,0,0,0),
 10: (0,0,0,0,-1,0,0),
 11: (0,0,0,-1,0,0,0),
 12: (-6,-4,-3,-3,-1,0,0),
 13: (-3,-2,-2,-1,-1,0,0),
 14: (3,2,2,1,1,0,0),
 15: (3,2,2,0,1,0,0)}
\end{verbatim}
\end{tiny}
\end{minipage}
\begin{minipage}{2.7 in}
The direction vectors of the tropical rays.  Since the 
homogeneity space is positive-dimensional, the directions are not uniquely specified. 
For instance, the vectors $(-5,-4,-3,-2,-1,0,0)$ and $(0,0,0,0,0,0,-1)$
represent the same ray. Note that \texttt{Gfan} uses negated weight vectors.
\end{minipage}

\vspace{.2cm}

\begin{minipage}{2 in}
\begin{tiny}    
\begin{verbatim}
Rays incident to each 
dimension 2 cone:
{{2,6}, {3,7},                                                   
{2,4}, {3,5},                                                   
{4,9}, {5,10},                                                   
{4,8}, {5,8},                                                  
{8,11},                                                    
{0,12}, {1,12},
{0,1},                                                
{1,6}, {0,7},                                                               
{1,9}, {0,10},                                                    
{0,13}, {1,13},                                                  
{6,14}, {7,14},                                                
{9,13}, {10,13},                                                
{6,10}, {7,9},                                             
{6,7},                                                  
{11,12},                                                               
{11,15},
{14,15}}
\end{verbatim}
\end{tiny}
\end{minipage}
\begin{minipage}{2.7 in}
  The cones in $\trop(I)$ are listed from highest to lowest dimension.
  Each cone is named by the set of rays on it.  There are 28 two-dimensional cones,
  broken down into 11 orbits of size 2 and 6 orbits of size~1.
\end{minipage}

\medskip
 
\noindent The further output, which is not displayed here,
shows that the 16 rays break down into 5 orbits of size 2 and 6
orbits of size 1.  
\begin{figure}
\epsfig{file=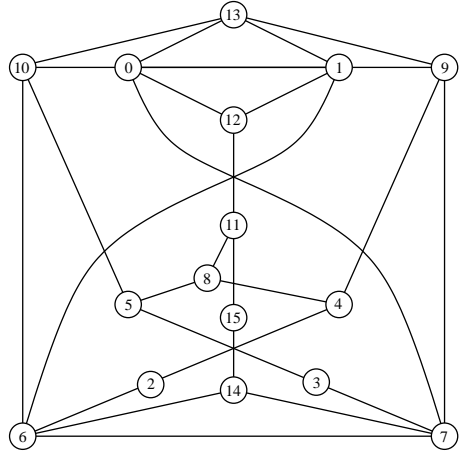,width=6cm}
\caption{The tropical variety of the ideal generated by the 
$3 \times 3$ minors of the generic $4 \times 4$ Hankel matrix.}
\label{fig: hankel}
\end{figure}
\end{example}  

Using the same procedure, we now compute several more examples.
\begin{example} \label{3minors of Hankel5by5}
Let $I$ be the ideal generated by the $3 \times 3$ minors of the
generic $5 \times 5$ Hankel matrix.  We again use the symmetry group
$\zz/2$.  The tropical variety is a graph with vertex degrees ranging
from 2 to 7.  
\begin{tiny}
\begin{verbatim}
Ambient dimension: 9
Dimension of homogeneity space: 2 
Dimension of tropical variety: 4
Simplicial: true
F-vector: (28,53) 
\end{verbatim}
\end{tiny}
\end{example}

\begin{example} \label{3minors of 3by5}
Let $I$ be the ideal generated by the $3 \times 3$ minors of a generic
$3 \times 5$ matrix.  We use the symmetry group $S_5 \times S_3$, where
$S_5$ acts by permuting the columns and $S_3$ by permuting the rows.  
\begin{tiny}
\begin{verbatim}
Ambient dimension: 15 
Dimension of homogeneity space: 7  
Dimension of tropical variety: 12 
Simplicial: true
F-vector: (45,315,930,1260,630)
\end{verbatim}
\end{tiny}
\end{example}

\begin{example} \label{3minors of nbynsymmetric}
  Let $I$ be the ideal generated by the $3 \times 3$ minors of a
  generic $4 \times 4$ symmetric matrix. We use the symmetry group
  $S_4$ which acts by simultaneously permuting the rows and the
  columns.
\begin{tiny}
\begin{verbatim}
Ambient dimension: 10 
Dimension of homogeneity space: 4 
Dimension of tropical variety: 7
Simplicial: true
F-vector: (20,75,75)
\end{verbatim}
\end{tiny}
If we take the $3 \times 3$ minors of a 
generic $5 \times 5$ symmetric matrix then we get
\begin{tiny}
\begin{verbatim}
Ambient dimension: 15 
Dimension of homogeneity space: 5 
Dimension of tropical variety: 9 
Simplicial: true
F-vector:    (75, 495, 1155, 855)
\end{verbatim}
\end{tiny}
\end{example}

\begin{example} \label{2by2 commuting}
Let $I$ be the prime ideal of a 
pair of commuting $2 \times 2$ matrices.  That is,
$I \subset \cc[a,b,\ldots,h] $ is defined by the matrix equation
$$
\left( \begin{array}{cc}
a & c \\
b & d 
\end{array} \right) 
\left( \begin{array}{cc}
e & g \\
f & h
\end{array} \right)
-
\left( \begin{array}{cc}
e & g \\
f & h
\end{array} \right)
\left( \begin{array}{cc}
a & c \\
b & d 
\end{array} \right) 
= 0.$$
 The tropical variety is the graph $K_4$,
 which \texttt{Gfan} reports as follows:
\begin{tiny}
\begin{verbatim}
Ambient dimension: 8
Dimension of homeogeneity space: 4
Dimension of tropical variety: 6
Simplicial: true
F-vector: (4,6)
\end{verbatim}
\end{tiny}
\vskip .1cm
If $I$ is the ideal of  $3 \times 3$ commuting symmetric matrices
then we get:
\begin{tiny}
\begin{verbatim}
Ambient dimension: 12
Dimension of homeogeneity space: 2
Dimension of tropical variety: 9
Simplicial: false
F-vector: (66,705,3246,7932,10888,8184,2745)
\end{verbatim}
\end{tiny}   
\end{example}


\section{Tropical variety versus Gr\"obner fan} \label{sec:gfan vs tropical variety}

In this paper we developed tools for computing the
 tropical variety $\TT(I)$ of a
 $d$-dimensional homogeneous prime ideal $I$
 in a polynomial ring $\cc[\xx]$. We took advantage of the
fact that, since $I$ is homogeneous, the set $\TT(I)$
has naturally the structure of a polyhedral fan, namely,
$\TT(I)$ is the collection of all cones in the Gr\"obner fan
of $I$ whose corresponding initial ideal is monomial-free.
A naive algorithm would be to compute the Gr\"obner fan of
$I$ and then retain only those $d$-dimensional cones who survive
the monomial test (Algorithm \ref{alg: monomial}).
The software {\tt Gfan} also computes the full Gr\"obner fan of $I$,
and so we tested this naive algorithm. We found it to be
too inefficient. The reason is that the vast 
majority of $d$-dimensional cones in the Gr\"obner fan
of $I$ are typically not in the tropical variety $\TT(I)$.

\begin{example}
Consider the ideal $I$ in Example \ref{ex:hankel4}
which is generated by the $3 \times 3$-minors of a generic
$4 \times 4$-Hankel matrix. Let $J = {\rm in}_w(I)$ be its initial ideal with respect to the first vector $w$ in the list of rays. The initial ideal $J$ defines a tropical curve consisting of five rays and the origin. The curve is a subfan of the much more complicated Gr\"obner fan of $J$. The Gr\"obner fan is full-dimensional in $\rr^7$ with $C_0(J)$ being three-dimensional. Its f-vector equals $(1,7167,32656,45072,19583)$. Of the $7167$ rays only $5$ are in the tropical variety. The Gr\"obner fan of $J$ is the link of the Gr\"obner fan of $I$ at $w$. We were unable to compute the full Gr\"obner fan of $I$. 
\end{example}

\begin{example} {\bf Toric Ideals.} \label{TI}
  Let $I= \langle \xx^\u - \xx^\v \,:\, A\u = A\v \rangle$ be the toric ideal of a matrix 
  $A \in \mathbb Z^{d \times n}$ of rank $d$. The ideal $I$ is a prime of
  dimension $d$. The tropical variety $\TT(I)$ coincides with the
   homogeneity space $C_0(I)$ which is just the row space of $A$.
 Hence $\TT(I)$ modulo $C_0(I)$ is a single point.
   Yet, the   Gr\"obner fan of $I$ can be very complicated, as it encodes the
  sensitivity information for an infinite family of integer programs  \cite[Chapter 7]{GBCP}. 
\end{example}

We next exhibit a family of ideals such that the number of
rays in $\TT(I)$ is constant while the
number of rays in the Gr\"obner fan of $I$ grows linearly.

\begin{theorem}
Fix $n \!=\! 3, d \!=\! 1$ and for any integer $p \geq 1$ consider the ideal 
$$ I_p \,\,\, = \,\,\, \langle \,\underline{x} - (z+1)^{p+2}, \, \underline{y} - (z-1)^p \,\rangle .$$
Then $\TT(I_p)$ consists of $4$ rays but
the Gr\"obner fan of $I_p$ has  $ \geq \frac{1}{4}(p + 1)$ rays.
\end{theorem}

\noindent {\sl Sketch of proof:}
The ideal $I_p$ is prime. Its variety is the parametric curve
$\, z \mapsto \bigl((z+1)^{p+2},(z-1)^p,z\bigr)$. The poles
and zeros of this map are $0,-1,+1,\infty $.
The tropical variety of  $I_p$ consists of the four rays defined by
the valuations at these points. These rays are generated by the columns of  
$$ \left( \begin{array}{cccc}
0 & 0 & p+2 & -p-2 \\
0 & p & 0   & -p   \\
1 & 0 & 0   & -1    
\end{array} \right).$$
We examine the Gr\"obner fan around the ray $\,w = -(p+2,p,1)$. The
initial ideal ${\rm in}_w(I_p)$ equals the toric ideal $\,\langle x -
z^{p+2}, y - z^p \rangle $.  To see this, we note that the two
generators of $I_p$ form a Gr\"obner basis with respect to the
underlined leading terms and ${\rm in}_w(I_p)$ is generated by ${\rm
  in}_w(g)$ for each $g$ in this Gr\"obner basis since $w$ lies in
this Gr\"obner cone. The Gr\"obner fan of ${\rm in}_w(I_p)$ is the
link at $w$ of the Gr\"obner fan of $I_p$. To prove the theorem we show
that the Gr\"obner fan of ${\rm in}_w(I_p)$ has at least
$\frac{1}{2}(p+1)$ distinct Gr\"obner cones. This implies, by Euler's
formula, that the Gr\"obner fan of ${\rm in}_w(I_p)$ has at least $
\frac{1}{4}(p + 1)$ rays and hence so does the Gr\"obner fan of $I_p$.

To argue that the Gr\"obner fan of ${\rm in}_w(I_p)$ has at least
$\frac{1}{2}(p+1)$ distinct Gr\"obner cones we use the methods in
\cite{GBCP}. More specifically, this involves first showing that the
binomials $g_j := y^j - z^{p-2(j-1)}x^{j-1}$ for $j = 1, \ldots,
\frac{p+1}{2}$ are all in the universal Gr\"obner basis of ${\rm
  in}_w(I_p)$. Each monomial in a binomial in the universal Gr\"obner
basis of a toric ideal contributes a minimal generator to some initial
ideal of the toric ideal. Thus there exist reduced Gr\"obner bases of
${\rm in}_w(I_p)$ in which the binomials $g_j$ are elements with
leading term $y^j$ for $j = 1, \ldots, \frac{p+1}{2}$. This implies
that these reduced Gr\"obner bases are all distinct, which completes
the proof.

 \qed

\smallskip

While the Gr\"obner fan is a fundamental object which has had a range
of applications (the Gr\"obner walk \cite{collart}, integer
programming (Example \ref{TI})), many computer algebra experts do not like it.
Their view is that {\em the Gr\"obner
  fan is a combinatorial artifact which is marginal to the real goal
  of computing the variety of $I$}.  While this opinion has some
merit, the story is entirely different for the subfan $\TT(I)$ of the
Gr\"obner fan.  In our view, {\em the tropical variety \underbar{is}
  the variety of $I$}.  Every point on $\TT(I)$ furnishes the starting
system for a {\em numerical homotopy} towards the complex variety of
$I$, see \cite[Chapter 3]{cbms}.  Thus computing $\TT(I)$ is not only
much more efficient than computing the Gr\"obner fan of $I$, it is
also geometrically more meaningful.

\bigskip

\noindent {\bf Acknowledgments.}
We thank Komei Fukuda for answering polyhedral computation questions
and for customizing {\tt cddlib} for our purpose, and we thank ETH
Z\"urich for hosting Anders Jensen and Bernd Sturmfels during the
summer of 2005. We acknowledge partial financial support by the
University of Aarhus, the Danish Research Training Council
(Forskeruddannelsesr\aa det, FUR), the Institute for Operations
Research at ETH, the Swiss National Science Foundation Project
200021-105202, and the US National Science Foundation (DMS-0456960,
DMS-0401047 and DMS-0354131).

\raggedright

\bigskip
\bigskip
\bigskip
 
\noindent Tristram Bogart, Department of Mathematics,
University of Washington, \break Seattle, WA 98195-4350, USA,
{\tt bogart@math.washington.edu}.
 
\bigskip
 
\noindent Anders Jensen, Institut for Matematiske Fag,
Aarhus Universitet, \break DK-8000 \AA rhus, Denmark,
{\tt ajensen@imf.au.dk}.
 
\bigskip
 
\noindent David Speyer, Department of Mathematics, University of Michigan \break Ann Arbor, MI 48109-1043, USA, {\tt speyer@umich.edu}.
 
\bigskip
 
\noindent Bernd Sturmfels, Department of Mathematics,
University of California, \break Berkeley, CA 94720-3840, USA,
{\tt bernd@math.berkeley.edu}.
 
\bigskip
 
\noindent Rekha Thomas, Department of Mathematics,
University of Washington, \break Seattle, WA 98195-4350, USA,
{\tt thomas@math.washington.edu}.

\end{document}